\documentclass[12pt, leqno]{article}
\usepackage[pdftex]{color, graphicx}
\usepackage{setspace}
\usepackage[utf8]{inputenc}
\usepackage{amsmath}
\usepackage[pdftex,bookmarks,colorlinks]{hyperref}
\usepackage{latexsym}
\usepackage{amsfonts}
\usepackage{amssymb}
\usepackage{amsthm}
\usepackage{paralist}
\usepackage{mathrsfs}
\usepackage{enumitem}
\usepackage[T1]{fontenc}
\usepackage{amsmath,amsthm}
\begin{document}
\newtheorem{thm}{Theorem}
\newtheorem{cor}[thm]{Corollary}
\newtheorem{prop}[thm]{Proposition}
\newtheorem{lem}{Lemma}
\theoremstyle{remark}\newtheorem{rem}{Remark}
\theoremstyle{definition}\newtheorem{defn}{Definition}


\title{Unconditional convergence of the differences of Fej\'er kernels on $L^2(\mathbb{R})$}

\author{Sakin Demir\\
Agri Ibrahim Cecen University\\ 
Faculty of Education\\
Department of Basic Education\\
04100 A\u{g}r{\i}, Turkey\\
E-mail: sakin.demir@gmail.com
}


\maketitle


\renewcommand{\thefootnote}{}



\renewcommand{\thefootnote}{\arabic{footnote}}
\setcounter{footnote}{0}


\begin{abstract}
Let $K_n(x)$ denote the Fej\'er kernel given by
$$K_n(x)=\sum_{j=-n}^n\left(1-\frac{|j|}{n+1}\right)e^{-ijx}$$
and  let $\sigma_nf(x)=(K_n\ast f)(x)$,  where as usual $f\ast g$ denotes the convolution of $f$ and $g$.\\
Let the sequence $\{n_k\}$ be lacunary. Then the series 
$$\mathcal{G}f(x)=\sum_{k=1}^\infty \left(\sigma_{n_{k+1}}f(x)-\sigma_{n_k}f(x)\right)$$ 
converges unconditionally for all $f\in L^2(\mathbb{R})$.\\
Let $(n_k)$ be a lacunary sequence, and $\{c_k\}_{k=1}^\infty \in \ell^\infty$. Define
$$\mathcal{R}f(x)=\sum_{k=1}^\infty c_k\left(\sigma_{n_{k+1}}f(x)-\sigma_{n_k}f(x)\right).$$ 
Then  there exists a constant $C>0$ such that
$$\|\mathcal{R}f\|_2\leq C\|f\|_2$$
for all $f\in L^2(\mathbb{R})$, i.e., $\mathcal{R}f$ is of strong type $(2,2)$. As a special case it follows that $\mathcal{G}f$ also is of strong type $(2,2)$.
\end{abstract}
\section{Preliminaries}
Even though the Fej\'er kernel has a long history in Fourier analysis, it is not hard to see by a quick literature review that this subject has not been studied extensively. For example, variation inequalities for the Fej\'er kernel have been studied in 2004 by R. L. Jones and G. Wang~\cite{rljgw}. Since then we do not see any remarkable work on this subject.  In this research we study the unconditional convergence of the  the Fej\'er kernel, we prove that the difference of the convolution with the  Fej\'er kernels for lacunary sequence converges unconditionally for all $f\in L^2(\mathbb{R})$. In order to prove our result we first control the Fourier transform and then use this control to prove required inequality for unconditional convergence.\\

\begin{defn} The series  $\sum_{n=1}^\infty x_n$  in a Banach space $X$ is said to converge unconditionally if the series  $\sum_{n=1}^\infty\epsilon_nx_n$  converges for all $\epsilon_n$ with $\epsilon_n=\pm 1$ for $n=1,2,3,\dots$.\\
The series  $\sum_{n=1}^\infty x_n$  in a Banach space $X$ is said to be weakly unconditionally convergent if for every functional $x^\ast\in X^\ast$ the scalar series $\sum_{n=1}^\infty x^\ast ( x_n)$  is unconditionally convergent.
\end{defn}
\begin{prop}\label{wuc} For a series $\sum_{n=1}^\infty x_n$  in a Banach space $X$ the following conditions are equivalent:
\begin{enumerate}[label=\}upshape(\roman*), leftmargin=*, widest=iii]
\item  The series  $\sum_{n=1}^\infty x_n$ is weakly unconditionally convergent;
\item There exists a constant $C$ such that for every $\{c_n\}_{n=1}^\infty\in \ell^\infty$
$$\sup_N\left\|\sum_{n=1}^N c_n x_n\right\|\leq C\|\{c_n\}\|_{\infty}.$$
\end{enumerate}
\end{prop}
\begin{proof}See page 59 in  P.~Wojtaszczyk~\cite{pwoj}.
\end{proof}
\begin{cor}\label{ucc}Let $X$ be a Banach space. If $\sum_{n=1}^\infty f_n$ is a series in $L^p(X)$, $1<p<\infty$, the following are equivalent:
\begin{enumerate}[label=\upshape(\roman*), leftmargin=*, widest=iii]
\item  The series  $\sum_{n=1}^\infty f_n$ is  unconditionally convergent;
\item There exists a constant $C$ such that for every $\{c_n\}_{n=1}^\infty\in \ell^\infty$
$$\sup_N\left\|\sum_{n=1}^N c_n f_n\right\|_p\leq C\|\{c_n\}\|_{\infty}.$$
\end{enumerate}
\end{cor}
\begin{proof}It is known (see page 66 in  P.~Wojtaszczyk~\cite{pwoj}) that every weakly unconditionally convergent series in a weakly sequentially complete space is unconditionally convergent. Since $L^p(X)$ is a  weakly sequentially complete space for $1<p<\infty$, the corollary follows from Proposition~\ref{wuc}.
\end{proof}
\begin{defn}\label{lacunary} A sequence $(n_k)$ of integers is called lacunary if there is a constant $\alpha >1$ such that
$$\frac{n_{k+1}}{n_k}\geq\alpha$$
for all $k=1,2,3,\dots$ .
\end{defn}
\section{The Results}
We denote by $K_n(x)$ the Fej\'er kernel given by
$$K_n(x)=\sum_{j=-n}^n\left(1-\frac{|j|}{n+1}\right)e^{-ijx}.$$

We let $\sigma_nf(x)=(K_n\ast f)(x)$,  where as usual $f\ast g$ denotes the convolution of $f$ and $g$.\\
\noindent
Our first result is the following:
\begin{thm}\label{occfk}Let the sequence $\{n_k\}$ be lacunary. Then the series 
$$\mathcal{G}f(x)=\sum_{k=1}^\infty \left(\sigma_{n_{k+1}}f(x)-\sigma_{n_k}f(x)\right)$$ 
converges unconditionally for all $f\in L^2(\mathbb{R})$.
\end{thm}

\begin{proof}
Let $\{c_k\}_{k=1}^\infty \in \ell^\infty$ and define
$$T_Nf(x)=\sum_{k=1}^Nc_k \left(\sigma_{n_{k+1}}f(x)-\sigma_{n_k}f(x)\right).$$
In order to prove that  $\mathcal{G}f$ converges unconditionally for all $f\in L^2(\mathbb{R})$ we have to show that for every  $\{c_n\}_{n=1}^\infty\in \ell^\infty$  there exists a constant $C>0$ such that
$$\sup_N\|T_Nf\|_2\leq C\|\{c_n\}\|_{\infty}$$
for all $f\in L^2(\mathbb{R})$ since this will verify the condition of  Corollary~\ref{ucc} for $\mathcal{G}f$.\\

Let 
$$S_N(x)=\sum_{k=1}^N\left(K_{n_{k+1}}(x)-K_{n_k}(x)\right).$$
\noindent
We clearly have
\begin{align*}
|\widehat{S}_N(x)|&=\left|\sum_{k=1}^N\left(\widehat{K}_{n_{k+1}}(x)-\widehat{K}_{n_k}(x)\right)\right|\\
&\leq \sum_{k=1}^N\left|\widehat{K}_{n_{k+1}}(x)-\widehat{K}_{n_k}(x)\right|.
\end{align*}
We first want to show that there exits a constant $C>0$ such that
$$ |\widehat{S}_N(x)|\leq C$$ 
for all $x\in\mathbb{R}$.\\
The Fej\'er kernel has a Fourier transform given by
$$
\widehat{K}_n(x) = \left\{ \begin{array}{ll}
1-\frac{|x|}{n+1}&\;\;\textrm{if}\:\:|x|\leq n;\\
0&\;\;\textrm{if}\:\:|x|>n.
\end{array} \right.
$$
Fix $x\in\mathbb{R}$, and let  $k_0$ be the first $k$ such that $|x|\leq n_k$ and let
$$I(x)=\sum_{k=1}^N\left|\widehat{K}_{n_{k+1}}(x)-\widehat{K}_{n_k}(x)\right|.$$
Then we have 
\begin{align*}
I(x)&=\sum_{k=1}^{n_{k_0-1}}\left|\widehat{K}_{n_{k+1}}(x)-\widehat{K}_{n_k}(x)\right|+\sum_{k=n_{k_0}}^N\left|\widehat{K}_{n_{k+1}}(x)-\widehat{K}_{n_k}(x)\right|\\
&=I_1(x)+I_2(x).
\end{align*}
\noindent
We clearly have $I_1(x)=0$ since $\widehat{K}_n(x)=0$ for $|x|>n$    so in order to control $|\widehat{S}_N(x)|$ it suffices to control
$$I_2(x)=\sum_{k=n_{k_0}}^N\left|\widehat{K}_{n_{k+1}}(x)-\widehat{K}_{n_k}(x)\right|.$$
We have
\begin{align*}
I_2(x)&=\sum_{k=n_{k_0}}^N\left|\widehat{K}_{n_{k+1}}(x)-\widehat{K}_{n_k}(x)\right|\\
&=\sum_{k=n_{k_0}}^N\left|1-\frac{|x|}{n_{k+1}+1}+\frac{|x|}{n_k+1}-1\right|\\
&=\sum_{k=n_{k_0}}^N\left|-\frac{|x|}{n_{k+1}+1}+\frac{|x|}{n_k+1}\right|\\
&\leq \sum_{k=n_{k_0}}^N\frac{|x|}{n_{k+1}+1}+\sum_{k=n_{k_0}}^N\frac{|x|}{n_k+1}\\
&\leq \sum_{k=n_{k_0}}^N\frac{|x|}{n_{k+1}}+\sum_{k=n_{k_0}}^N\frac{|x|}{n_k}\\
&\leq \sum_{k=n_{k_0}}^N\frac{n_{k_0}}{n_{k+1}}+\sum_{k=n_{k_0}}^N\frac{n_{k_0}}{n_k}.
\end{align*}
On the other hand, since the sequence $\{n_k\}$ is lacunary there is a real number $\alpha >1$ such that
$$\frac{n_{k+1}}{n_k}\geq \alpha$$
for all $k\in\mathbb{N}$. Hence we have
$$\frac{n_{k_0}}{n_k}=\frac{n_{k_0}}{n_{k_0+1}}\cdot\frac{n_{k_0+1}}{n_{k_0+2}}\cdot\frac{n_{k_0+2}}{n_{k_0+3}}\cdots \frac{n_{k-1}}{n_k}\leq\frac{1}{\alpha^k}.$$
Thus we get
$$\sum_{k=n_{k_0}}^N\frac{n_{k_0}}{n_k}\leq \sum_{k=n_{k_0}}^N\frac{1}{\alpha^k}\leq \frac{\alpha}{\alpha -1}.$$
and similarly, we have
$$\sum_{k=n_{k_0}}^N\frac{n_{k_0}}{n_{k+1}}\leq \frac{\alpha}{\alpha -1}$$
and this proves that
$$I_2(x)\leq 2\frac{\alpha}{\alpha -1}.$$
Since the bound does not depend on the choice of $x\in \mathbb{R}$ what we have just proved is true for all $x\in \mathbb{R}$.\\
We conclude that there exits a constant $C>0$ such that
$$|\widehat{S}_N(x)|\leq C\;\;\;\;\;\;\;\;\;\;\;\;\; \textrm{($\ast$)}$$
for all $x\in\mathbb{R}$ and $N\in\mathbb{N}$.\\
We now have
\begin{align*}
\|T_Nf\|_2^2&=\int_{\mathbb{R}}\left|\sum_{k=1}^Nc_k\left(\sigma_{n_{k+1}}f(x)-\sigma_{n_k}f(x)\right)\right|^2\, dx\\
&=\int_{\mathbb{R}}\left|\sum_{k=1}^Nc_k\left(K_{n_{k+1}}\ast f(x)-K_{n_k}\ast f(x)\right)\right|^2\, dx\\
&\leq \|\{c_n\}\|_{\infty}^2 \int_{\mathbb{R}}\left|\sum_{k=1}^N\left(K_{n_{k+1}}\ast f(x)-K_{n_k}\ast f(x)\right)\right|^2\, dx\\
&= \|\{c_n\}\|_{\infty}^2 \int_{\mathbb{R}}|S_N\ast f(x)|^2\, dx\\
&= \|\{c_n\}\|_{\infty}^2 \int_{\mathbb{R}}|\widehat{S_N\ast f}(x)|^2\, dx\;\;\;\textrm{(by Plancherel's theorem)}\\
&= \|\{c_n\}\|_{\infty}^2\int_{\mathbb{R}}|\widehat{S}_N(x)|^2\cdot|\hat{ f}(x)|^2\, dx\\
&\leq C \|\{c_n\}\|_{\infty}^2\int_{\mathbb{R}}|\hat{ f}(x)|^2\, dx\;\;\;\textrm{(by ($\ast$))}\\
&=C\|\{c_n\}\|_{\infty}^2\int_{\mathbb{R}}| f(x)|^2\, dx\;\;\;\textrm{(by Plancherel's theorem)}\\
&=C\|\{c_n\}\|_{\infty}^2\|f\|_2^2
\end{align*}
and thus we get
$$\sup_N\|T_Nf\|_2\leq \sqrt{C}\|\{c_n\}\|_{\infty}\|f\|_2$$
which completes our proof.
\end{proof}
\begin{thm}\label{rfl2} Let $(n_k)$ be a lacunary sequence, and $\{c_k\}_{k=1}^\infty \in \ell^\infty$. Define
$$\mathcal{R}f(x)=\sum_{k=1}^\infty c_k\left(\sigma_{n_{k+1}}f(x)-\sigma_{n_k}f(x)\right).$$ 
Then  there exists a constant $C>0$ such that
$$\|\mathcal{R}f\|_2\leq C\|f\|_2$$
for all $f\in L^2(\mathbb{R})$, i.e., $\mathcal{R}f$ is of strong type $(2,2)$.
\end{thm}
\begin{proof}We have proved that in the proof of Theorem~\ref{occfk} that  given $N\in\mathbb{N}$ there exists a constant $C_1>0$ such that
$$\sum_{k=1}^N\left|\widehat{K}_{n_{k+1}}(x)-\widehat{K}_{n_k}(x)\right|\leq C_1$$
for all $x\in\mathbb{R}$, we also have by taking limit
\begin{align*}
\sum_{k=1}^\infty\left|\widehat{K}_{n_{k+1}}(x)-\widehat{K}_{n_k}(x)\right|^2&\leq \sum_{k=1}^\infty\left|\widehat{K}_{n_{k+1}}(x)-\widehat{K}_{n_k}(x)\right|\\
&\leq C_1
\end{align*}
$x\in\mathbb{R}$.\\
Then we obtain
\begin{align*}
\|\mathcal{R}f\|_2^2&=\int_{\mathbb{R}}\left|\sum_{k=1}^Nc_k\left(\sigma_{n_{k+1}}f(x)-\sigma_{n_k}f(x)\right)\right|^2\, dx\\
&=\int_{\mathbb{R}}\left|\sum_{k=1}^{\infty}c_k\left(K_{n_{k+1}}\ast f(x)-K_{n_k}\ast f(x)\right)\right|^2\, dx\\
&\leq \|\{c_n\}\|_{\infty}^2 \int_{\mathbb{R}}\left|\sum_{k=1}^\infty\left(K_{n_{k+1}}\ast f(x)-K_{n_k}\ast f(x)\right)\right|^2\, dx\\
&= \|\{c_n\}\|_{\infty}^2 \int_{\mathbb{R}}\sum_{k=1}^\infty\left|\left(K_{n_{k+1}}\ast f(x)-K_{n_k}\ast f(x)\right)\right|^2\, dx\\
&= \|\{c_n\}\|_{\infty}^2 \sum_{k=1}^\infty\int_{\mathbb{R}}\left|\left(K_{n_{k+1}}\ast f(x)-K_{n_k}\ast f(x)\right)\right|^2\, dx\\
&= \|\{c_n\}\|_{\infty}^2 \sum_{k=1}^\infty\int_{\mathbb{R}}\left|\left(\widehat{K_{n_{k+1}}\ast f}(x)-\widehat{K_{n_k}\ast f}(x)\right)\right|^2\, dx\;\;\;\textrm{(by Plancherel's theorem)}\\
&= \|\{c_n\}\|_{\infty}^2 \sum_{k=1}^\infty\int_{\mathbb{R}}\left|\widehat{K}_{n_{k+1}}(x)-\widehat{K}_{n_k}(x)\right|^2|\hat{f}(x)|^2\, dx\\
&= \|\{c_n\}\|_{\infty}^2\int_{\mathbb{R}} \sum_{k=1}^\infty\left|\widehat{K}_{n_{k+1}}(x)-\widehat{K}_{n_k}(x)\right|^2|\hat{f}(x)|^2\, dx\\
&\leq \|\{c_n\}\|_{\infty}^2C_1\int_{\mathbb{R}} |\hat{f}(x)|^2\, dx\\
&=\|\{c_n\}\|_{\infty}^2C_1\int_{\mathbb{R}}| f(x)|^2\, dx\;\;\;\textrm{(by Plancherel's theorem)}\\
&=\|\{c_n\}\|_{\infty}^2C_1\|f\|_2^2.
\end{align*}
This means that  there exists a constant $C>0$ such that
$$\|\mathcal{R}f\|_2\leq C\|f\|_2$$
for all $f\in L^2(\mathbb{R})$, i.e., $\mathcal{R}f$ is of strong type $(2,2)$.
\end{proof}
\begin{cor}Let $(n_k)$ be a lacunary sequence. Then  there exists a constant $C>0$ such that
$$\|\mathcal{R}f\|_2\leq C\|f\|_2$$
for all $f\in L^2(\mathbb{R})$, i.e., $\mathcal{R}f$ is of strong type $(2,2)$.
\end{cor}
\begin{proof}When we choose $c_k=1$ for all $k$ in the definition $\mathcal{R}f$ we obtain 
$$\mathcal{G}f=\mathcal{R}f$$
and the proof follows from Theorem~\ref{rfl2}.
\end{proof}

\end{document}